\documentstyle[11pt]{article}
\def\hybrid{\topmargin 0pt      \oddsidemargin 0pt
        \headheight 0pt \headsep 0pt
        \textwidth 160true mm       
        \textheight 231true mm         
        \marginparwidth 0.0in
        \parskip 0pt plus 1pt   \jot = 1.5ex}
\catcode`\@=11
\def\marginnote#1{}

\newcount\hour
\newcount\minute
\newtoks\amorpm
\hour=\time\divide\hour by60
\minute=\time{\multiply\hour by60 \global\advance\minute by-\hour}
\edef\standardtime{{\ifnum\hour<12 \global\amorpm={am}%
        \else\global\amorpm={pm}\advance\hour by-12 \fi
        \ifnum\hour=0 \hour=12 \fi
        \number\hour:\ifnum\minute<10 0\fi\number\minute\the\amorpm}}
\edef\militarytime{\number\hour:\ifnum\minute<10 0\fi\number\minute}

\def\draftlabel#1{{\@bsphack\if@filesw {\let\thepage\relax
   \xdef\@gtempa{\write\@auxout{\string
      \newlabel{#1}{{\@currentlabel}{\thepage}}}}}\@gtempa
   \if@nobreak \ifvmode\nobreak\fi\fi\fi\@esphack}
        \gdef\@eqnlabel{#1}}
\def\@eqnlabel{}
\def\@vacuum{}
\def\draftmarginnote#1{\marginpar{\raggedright\scriptsize\tt#1}}

\def\draft{\oddsidemargin -.5truein
        \def\@oddfoot{\sl preliminary draft \hfil
        \rm\thepage\hfil\sl\today\quad\militarytime}
        \let\@evenfoot\@oddfoot \overfullrule 3pt
        \let\label=\draftlabel
        \let\marginnote=\draftmarginnote
   \def\@eqnnum{(\theequation)\rlap{\kern\marginparsep\tt\@eqnlabel}%
\global\let\@eqnlabel\@vacuum}  }
\relax

\hybrid

\def\cp{{\bf CP}^n} 
\def\gggg{{\bf g}}
\title{Double quantization of  $\cp$ type orbits by generalized Verma modules}

\author{J.~Donin,\\
Department of Mathematics, Bar-Ilan University,\\
52900 Ramat-Gan, Israel,\\
D.~Gurevich,\\
ISTV, Universit\'e de Valenciennes\\
59304 Valenciennes, France\\
S.~Khoroshkin,\\
Institute of Theoretical and Experimental Physics\\
117259 Moscow, Russia}

\begin{document}

\maketitle
\begin{abstract}
It is known that symmetric orbits in ${\bf g}^*$ for any simple Lie algebra
${\bf g}$ are equiped with a Poisson pencil generated by 
 the Kirillov-Kostant-Souriau bracket and the reduced Sklyanin bracket
associated to the "canonical" R-matrix. We realize quantization of this 
Poisson 
pencil on $\cp$ type orbits (i.e. orbits in $sl(n+1)^*$ whose 
real compact form is ${\bf CP}^n$)
by means of q-deformed Verma modules.
\end{abstract}

{\bf AMS Mathematics Subject Classification, 1991 :} 17B37, 81R50

{\bf Key words :} Poisson bracket, Poisson pencil, R-matrix bracket, (double) 
quantization, 
flat deformation, orbit of ${\bf CP}^n$ type, (generalized) Verma module, 
quantum group,
(twisted) Hopf algebra, braided algebra, braided module

\def\Uq{U_q(sl(n))}
\def\ugqh{U(\gq)_{\h}}
\def\ugnh{U(\gggg_{\nu})_{\h}}
\def\iqma{I_q(\mu,\alpha)}
\def\aqma{{\aaaa}_q(\mu,\alpha)}
\def\ugg{U(\gggg)}
\def\ugn{U(\gggg)[[\nu]]}
\def\ug2{U(\gggg)^{\ot 2}}
\def\fn{F_{\nu}}
\def\oH{{\overline H}}
\def\bc{{\bf C}}
\def\alz{\alpha_0}
\def\cpp{{\bf CP}^{n-1}}
\def\oDe{\overline{\Delta}}
\def\De{\Delta}
\def\uq{U_q({\bf g})}
\def\GG{{\gggg}^{\ot 2}}
\def\gq{{\bf g}_q}
\def\gg2{{\bf g}^{\ot 2}}
\def\gqq{{\gggg}_q^{\ot 2}}
\def\cs{\cal S}
\def\csz{{\cal S}_0}
\def\Ob{Ob\,(\cs)}
\def\Ip{I_+}
\def\os{\overline s}
\def\Im{I_-}
\def\Ipm{I_{\pm}}
\def\Ipq{I^q_+}
\def\Imq{I^q_-}
\def\Ipmq{I^q_{\pm}}
\def\h{{\hbar}}
\def\la{\lambda}
\def\co{{\cal O}}
\def\coom{{\cal O}_{\omega}}
\def\a{\alpha}
\def\ch{{\bf C}[[\h]]}
\def\cn{{\bf C}[[\nu]]}
\def\hhhh{{\bf h}}
\def\kkkk{{\bf k}}
\def\nnnn{{\bf n}}
\def\pppp{{\bf p}}
\def\kkkk{{\bf k}}
\def\g{{\bf g}}
\def\vom{V_{\om}}
\def\orom{\overline{\rho}_{\om}}
\def\mom{M_{\om}}
\def\ve{\varepsilon}
\def\ep{\epsilon}
\def\bea{\begin{eqnarray}}
\def\eea{\end{eqnarray}}
\def\beq{\begin{equation}}          \def\bn{\beq}
\def\eeq{\end{equation}}            \def\ed{\eeq}
\def\nn{\nonumber}                 
\def\ee{{\rm e}}
\def\aaaa{\cal A}
\def\aaa{\cal A}
\def\ah{{\cal A}_{\h}}
\def\aq{{\cal A}_{q}}
\def\an{{\cal A}_{\nu}}
\def\ahn{{\cal A}_{\h,\nu}}
\def\ahq{{\cal A}_{\h,q}}
\def\ot{\otimes}
\def\om{\omega}
\def\Vma{V_{\mu,\alpha}}
\def\ad{{\rm ad}}                   \def\gij{g_{i,j}}
\def\gkl{g_{k,l}} 
\def\Tg{T(\gggg)}
\def\tS{\tilde{S}}
\def\Tgq{T(\gq)}
\def\gq{{\bf g}_q}
\def\Ml{M^{\lambda}}
\def\TTg{T_{\mu,\a}(\gggg)}
\def\Vm{V_{\mu\omega_1}}
\def\End{{\rm End\, }}
\def\Ker{{\rm Ker\, }}
\def\Im{{\rm Im\, }}
\def\id{{\rm id\, }} 
\def\Vect{{\rm Vect\,}}
\def\Fun{{\rm Fun\,}}
\def\vv{V^{\ot 2}}
\renewcommand{\theequation}{{\thesection}.{\arabic{equation}}}
\def\r#1{\mbox{(}\ref{#1}\mbox{)}}
\setcounter{equation}{0}
\newtheorem{proposition}{Proposition}
\newtheorem{corollary}{Corollary}
\newtheorem{theorem}{Theorem}
\newtheorem{definition}{Definition}
\newtheorem{remark}{Remark}

\def\R{{\cal R}}
\def\n{{\bf n}}   \def\va{\varepsilon}
\def\so{s_{\omega_1+\omega_{n-1}}}

\section{Introduction}

The problem of  quantization of 
Poisson brackets is one of the most important in  
mathematical physics. In the framework of the deformation quantization scheme
going back to the works by A. Lichneriwicz and his school (cf. \cite{BFFLS}) 
it can be
formulated as follows. Given a variety $M$ equipped with a Poisson bracket,
it is necesary to construct a flat deformation 
$\ah$ of an  algebra   ${\aaaa}={\Fun}{(M)}$ of functions over 
$M$\footnote{We say  that an algebra $\ah$ depending on a formal
parameter $\h$ is a {\em flat
deformation}  (or simply, deformation) of $\aaaa$ if 
${\aaa}={\ah}/\h {\ah}$  and $\ah$
is isomorphic to ${\aaaa}[[\h]]$
as $\ch$-modules. Hereafter $V[[\h]]$ where $V$ is a linear space stands 
for the completion in
$\h$-adic topology of 
$V\ot_{\bf C}\ch$ (in what follows the basic field is $k={\bf C}$).  
Abusing notation 
we will let ${\aaaa}\to{\ah}$ denote the deformation in question. 
Two parameter 
flat deformation  can be defined in a similar way.} such that the 
corresponding Poisson bracket (which exists for any flat deformation of 
a commutative algebra) coincides with the initial one.

The existence of such a quantization for any nowhere degenerated 
(i.e. defined by a 
symplectic structure) Poisson bracket had been shown in \cite{DL}. Recently, M.
Kontsevich \cite{K} has proved that any Poisson bracket is quantizable 
in the above sense.

Nevertheless, physysists are  interested in an operator quantization, i.e. they
want to realize the quantum algebra $\ah$  as an 
operator algebra in a linear  (ideally, Hilbert) space. This enables
 them to carry out
a spectral analysis of Hamiltonians and to compute partition functions and 
other
numerical characteristics of quantum models. Such a quantization of 
nondegenerated 
Poisson bracket (on any compact smooth variety) has been realized by B.Fedosov
\cite{F}. In fact, the famous Kirillov-Duflo orbit method which 
consists in assigning a representation $\rho:\gggg\to \End V$ of a Lie 
algebra $\gggg$ 
to an orbit $\co\subset\gggg^*$ can be considered as a particular case of 
the Fedosov approach.
(We do not discuss here the limits of the orbit method, in the sequel we 
will restrict ourselves to semisimple orbits in $\gggg^*$ for simple Lie 
algebras $\gggg$).

The quantization procedure suggested by Fedosov  leads to an operator 
algebra equipped with a
commutative trace. In fact, such a trace is delivered for appropriated 
quantum algebras by
the Liouville mesure of the initial Poisson bracket. However, a generic 
Poisson bracket
does not possess any invariant measure and consequently it is not clear 
what is a trace
in the corresponding quantum algebra.

In the earlier 90's one of the authors (D.G.) suggested certain Poisson 
brackets
associated to classical R-matrices whose quantization leads to operator
algebras in twisted categories. Traces in such algebras are also twisted
(cf. \cite{G3}, \cite{GRZ}). These 
algebras arise from quantization of Poisson pencils generated by the 
linear Poisson-Lie bracket on $\gggg^*$ or its restriction to an orbit,
called the Kirillov-Kostant-Souriau (KKS) one, 
and by a bracket naturally associated to
a solution $R\in \bigwedge^2(\gggg)$ 
of the classical {\em nonmodified} Yang-Baxter equation
\beq
[[R,R]]=[R^{12},\,R^{13}]+[R^{12},\,R^{23}]+[R^{13},\,R^{23}]=0.
\label{cyb}
\eeq
Here, as usual $R^{12}=R\ot \id$ etc.

Let us describe the latter bracket. Let $M$ be a variety equipped with a 
representation
$\rho:\gggg\to \Vect(M)$ where $\Vect (M)$ stands for the
space of vector fields. Then the following bracket
\beq
\{f,g\}_R=\mu\,<\rho^{\ot 2}(R), df\ot dg>,\,\,f,g \in \Fun(M)
\label{rmb}
\eeq
is Poisson.
Here  $<\,\,,\,\,>$ is the pairing between the differential forms and 
vector fields on
$M$  extended to their tenzor powers and $\mu$ is the usual
commutative product in the space $\Fun(M)$.
The bracket $\{\,\,,\,\,\}_R$ is called the {\em R-matrix bracket}.
If $M=\gggg^*$ or $M=\co\subset\gggg^*$ is an orbit we take as $\rho$ 
the coadjoint 
representation or its restriction to the orbit.

It is not difficult to see that in the latter case any bracket of the family
\begin{equation}
\{\,\,,\,\,\}_{a,b}=a\{\,\,,\,\,\}_{KKS}+b\{\,\,,\,\,\}_R \label{pp}
\end{equation}
is Poisson. Here by $\{\,\,,\,\,\}_{KKS}$ we mean either the KKS bracket or 
the linear  Poisson-Lie
one on $\gggg^*$.
Thus, we have a Poisson pencil well defined on $\gggg^*$ or on an orbit in 
$\gggg^*$.

A procedure of quantizing this Poisson pencil can be realized in two steps 
(we call 
such procedures "double quantization"). In the first step one quantizes only
the KKS bracket by means of the orbit method or by means of generalized Verma
modules as it is discribed in Section 3. Then one twists the quantum
operator algebra as it is described in Section 2. The resulting object is
a two parameter operator algebras  in a twisted category. It comes with 
a deformed trace which is no longer commutative but S-commutative 
in the spirit of a super-trace. Here $S$ is an involutive 
$(S^2=id)$  twist, i.e. an operator acting in tensor square of this algebra 
and satisfying  
quantum Yang-Baxter equation (QYBE)
$$S^{12}S^{23}S^{12}=S^{23}S^{12}S^{23}.$$

Let us remark that by means of a similar twisting one can introduce natural 
"S-analogues"  of basic objects of geometry and analysis. Thus, S-analogues of 
commutative algebras, vector fields, Lie algebras, (formal) Lie groups 
were defined,
in a spirit of super-theory,  in \cite{G1},  \cite{G4} (cf. also 
\cite{GRR}, \cite{GRZ}).\footnote{All these objects
are also well defined for some nonquasiclassical twists 
 (i.e.  twists which cannot be obtained by a deformation of the ordinary
flip $S=\sigma$). One can naturally associate to an involutive  twist
$S:\vv\to\vv$ where $V$ is a linear space $S$-symmetric and $S$-skewsymmetric
algebras of $V$. The Poincar\'e series of these algebras related to 
quasiclassical
twists, i.e. those being deformation of the flip,  coincide with classical ones
while those series corresponding to nonquasiclassical twists can
differ drastically from the classical (and super-) ones. The first
examples of such twists were given in \cite{G1} (similar twists of Hecke 
type has been
introduced in \cite{G2}).}. However the straightforward generalization
of these notions to noninvolutive twists (connected, say, to the 
quantum group (QG) $\uq$) leads
as a rule to a nonflat deformation. (From our viewpoint the principle "raison
d'\^etre" for objects
belogning to the category of $\uq$-modules is that they should represent a flat
 deformation of their classical counterparts.)

The main purpose of the paper is to realize an operator quantization of 
Poisson pencils
\r{pp} associated to the "canonical"  classical R-matrix
\beq
R=\sum_{\alpha\in\Omega^+}\frac{X_{\alpha}\wedge X_{-\alpha}}
{<X_{\alpha}, X_{-\alpha}>}\in{\bigwedge}^2(\gggg),
\label{rm}
\eeq 
where $\gggg$ is a complex simple Lie algebra,
$\Omega^+$ stands for the set of its positive roots with respect 
to a fixed triangular decomposition of $\gggg$ and
$<\,\,,\,\,>$ stands for the Killing form. 

This R-matrix satisfies the so-called classical {\em modified} YBE which means 
that the element
$[[R,R]]$ is nontrivial and $\gggg$-invariant. Since this element is not 
identically zero 
the associated R-matrix bracket is 
 Poisson only on  varieties where the  three-vector field
$\rho^{\ot 3}([[R,R]])$ vanishes.
Such  varieties were called in
\cite{GP}  {\em R-matrix type} ones. 
All R-matrix type orbits in $\gggg^*$  
 were
classified in \cite{GP}. In particular, all symmetric orbits in $\gggg^*$ are 
 R-matrix type  varieties.
(Let us recall that an orbit $\co_x$ of a point $x$ is called 
symmetric if there exists a  decomposition 
$\gggg=\kkkk\oplus{\bf m}$ 
where $\kkkk$ is the 
stabilizer of $x$ such that
$[\kkkk,\kkkk]\subset\kkkk,\,\, [{\bf m},{\bf m}]\subset\kkkk,\,\,[\kkkk, 
{\bf m}]\subset{\bf m}).
$

Moreover, the R-matrix bracket  over a symmetric orbit 
coincides with one of the two left- or right- invariant components  
of the  Sklyanin bracket reduced to the orbit (recall that the 
Sklyanin bracket  equals to  a difference
between  left- and right- invariant brackets defined by \r{rmb}, where
$\rho$ is the natural homomorphism of 
 $\gggg$ into  the space of left- or right- invariant vector fields 
on the corresponding group $G$).
 Meanwhile, the other  component being reduced becomes proportional to
the KKS bracket. This implies that on any symmetric orbit the R-matrix bracket
 and the KKS one are compatible and therefore they generate the Poisson pencil 
\r{pp}. 

Note that the one-sided invariant components of the Sklyanin bracket can be 
reduced to
any semisimple (ss) orbit in $\gggg^*$ (i.e., to that of an ss element),  but 
each of them becomes Poisson brackets only on symmetric orbits (cf. 
\cite{KRR} and \cite{DG2}).

The Poisson pencil \r{pp} 
with R-matrix \r{rm} on symmetric orbits has been
 quantized in the spirit of deformation quantization in \cite{DS}. 
 The resulting object
of the quantization procedure suggested in \cite{DS} is a two parameter family
 of associative  $\uq$-invariant algebras. Let us make precise that
an algebra $\aaaa$ is called $\uq$-invariant (or $\uq$-covariant) if
\beq
u\cdot (x_1x_2)=\left( u_{(1)}\cdot x_1\right)
\left( u_{(2)}\cdot x_2\right),\,\,\forall u \in \uq,\,\, x_1,x_2 \in \aaaa.
\label{1.1}
\eeq
Hereafter 
$u_{(1)}\ot u_{(2)}$ stands for $\De(u)$ (the Sweedler's notations).
Algebras $\aaaa$ 
possessing this property will be called {\em quantum or braided}
 ones, while by {\em twisted algebras} we mean the algebras beloning to a  
twisted category. 

However, the product in the quantum algebra is realized in \cite{DS} 
by a series in two 
formal parameters, meanwhile the QG $\uq$ appears    
as $U(\gggg)[[\nu]]$ but equipped with a deformed coproduct (the 
so-called, Drinfeld's realization,
cf. Section 2).

In the present paper we perform a double quantization of $\cp$ type orbits 
by an operator method. 
By $\cp$ type orbits we mean the orbits in $sl(n)^*$ of elements
$\mu\om_1$ or $\mu\om_{n-1}$ where $\om_1\,\,(\om_{n-1})$ is the first 
(the last) fundamental weight of $sl(n)$ and $\mu\in{\bf C}$ is an arbitary 
nontrivial factor.
 Compact forms of theses complex orbits are 
just ${\bf CP}^{n-1}$ embedded as closed algebraic varieties in $su(n)^*$.

More precisely, we represent our two parameter quantum object
$\ahq$ as operator algebra in braided (or q-deformed) generalized Verma 
modules. Similary to
 the previous case arising from the classical nonmodified YBE,
 our quantization procedure consists of two steps.

The first, "classical", step is realized as follows.
There exists a natural way to quantize ss orbits in $\gggg^*$ for any 
simple Lie algebra $\gggg$ by generalized Verma modules. 
Let $\mom$ be such a module of highest weight $\om$
(its construction  is given 
in Section 3) and $\rho_{\om} :T(\gggg)\to \End\,\mom$ be the corresponding 
representation of the free tensor algebra $T(\gggg)$.
Then,  the operator algebras
$\ah=\Im\,\h \rho_{\om/\h}\subset\End \mom[[\h]]$
 can be treated as a quantum object with respect to 
the KKS bracket on the orbit  $\coom\subset\gggg^*$ of the element $\om$ 
(we regard $\om$ as an element of $\gggg^*$ as explained below).
The passage from the representation $\rho_{\om}$ to that $\h \rho_{\om/\h}$ 
will be reffered
below as a renormalization procedure. 

Let us remark that operator algebra $\ah$ is an object of the category of 
$\gggg$-invariant algebras similarly to the initial function algebra
$\aaaa=\Fun(\coom)$.

The second step consists in a braiding of the algebras $\ah$. As a result we 
get 
the mentioned above two parameter $\uq$-invariant operator algebra $\ahq$.
Let us emphasize that our approach to represent quantum algebras by means of 
braided
generalized Verma modules has the following advantage. 
The parameters $\h$ and $q$ can be specialized: the operator realization 
of the algebra $\ahq$ is well defined for any value of $\h$ and a generic
$q$.

Moreover, the flatness of deformation
${\aaaa}\to\ahq$ is assured automatically.
Let us remark that the quantum algebra $\ahq$ can be also
represented by a system of some algebraic equations. For the $\cp$ type orbits
these equations are quadratic-linear-constant. It is not so difficult to 
guess a general 
form of these equations. The problem is 
to find the exact meaning of factors occurring in them
which ensure flatness of deformation of the corresponding quotient algebras.
Different ways to look for these factors were discussed in
\cite{DG2}, \cite{GR}, \cite{DGR2}, \cite{G5}.
The  "operator method" presented here is  the most adequate way to solve 
this problem.

Thus, compared with  \cite{DS}  our approach enables us 
to realize quantum
conterparts of  the Poisson pencil  in question explicitely in the spirit 
of noncommutative
algebraic geometry.

The paper is organized as follows. In the next section we describe different
algebraic structures connected to involutive twists  arising 
from quantization of R-matrices satisfying  the
classical {\em nonmodified} YBE. 
We  
show  that certain quotients of
twisted Hopf algebras  are the
appropriate objects to describe explicitely quantized orbits in $\gggg^*$. 
We also analyze  the difference between 
this case refered in what follows as triangular or involutive
and that connected to the quasitriangular QG $\uq$. 

Section 3 is devoted to the "classical step" of quantization. The final object
of this step is the mentioned above algebras $\ah$. Then we realize  a
q-deformation of these algebras as follows. We equip $\gggg=sl(n)$ with a
structure of a
 $U_q(\gggg)$-module, extend the action of the QG $U_q(\gggg)$ to its 
enveloping
algebra and represent the latter algebra in the q-deformed 
generalized Verma modules considered on the first step. These constructions
are discribed in Sections 4 and 5. They result in 
a two parameter family $\ahq$ presented in the last section.

Completing the  introduction we want to put the following question: 
 how is it possible 
to define a proper trace in a quantum algebra arising in virtue of \cite{K}
from a given Poisson bracket? As our examples show, such traces
are not necessary commutative. (Although we are dealing with the 
complexification 
of $\cp$ the trace 
defined by a projection of the algebra $\ahq$ to its trivial component 
is well defined on this algebra since it corresponds
to the compact form of the orbits in question,  cf.
\cite{GV} where such a trace in $sl(2)$ case is studied).

{\bf Acknowledgements}  The autors are deeply indebted to S. Majid who read
a preliminary version of the paper and made presious remarks. One  of the 
authors (S.K.) 
 was supported by grants RFBR 96-01-01421 and INTAS 93-10183.
 He would like also to thank Universit\'e Lille-1 (France) for the warm 
hospitality during
 his stay when this paper was started.

\setcounter{equation}{0}
\section{Triangular and quasitriangular cases: comparative description}

{\bf 2.1.} Let us first consider  certain algebraic structures  
arising from R-matrices  satisfying  the classical YBE \r{cyb}.
Let us fix such an R-matrix $R$.

By Drinfeld's result \cite{D1} there exists a series $F=\fn\in \ug2[[\nu]]$ 
quantizing the R-matrix $R$ in the following sense
$\fn=1+\nu P+...$ where $P\in \gggg^{\ot 2},\,\,P-P^{21}= 2R$, and
\beq
\De^{12}F\, F^{12}=\De^{23}F\, F^{23},\,\,\,\ve^1F=\ve^2F=1.
\label{def}
\eeq
Here $\De: U(\gggg)\to  U(\gggg)^{\ot 2}$ is the usual coproduct 
and 
$\ve: \ugg\to \bc$ is the counit in $\ugg$ (all operators are assumed to 
be extended  to $\ugn$ in a natural way). 
Using $F$ one can deform the usual Hopf structure of the algebra
$\ugg$ in (at least) two different ways. 

The first way consists of the following procedure. Let us introduce a new 
coproduct  setting
$$\De_F(u)=F^{-1}\De(u)F=F^{-1}_{(1)}u_{(1)}F_{(1)}\ot F^{-1}_{(2)}u_{(2)}
F_{(2)}.$$
Here $F_{(1)}\ot F_{(2)}$ (reps., $F_{(1)}^{-1}\ot F_{(2)}^{-1}$) stands
for $F$ (resp., $F^{-1}$) .

Then the algebra $\ugn$ equipped with the initial product and unit, the
coproduct $\De_F$ and the uniquely defined  counit and antipode 
(cf. \cite{D2} and
\cite{GM} where the antipode is expressed via $F$) becomes an Hopf 
algebra looking like the famous QG $\uq$. Let us denote it $H$.

Another way consists in simultaneous deformation of the product and coproduct
as follows
$$\overline{\De}(u)=ad\,F^{-1}(\De(u))=ad\,F_{(1)}^{-1}(u_{(1)})\ot 
ad\,F_{(2)}^{-1}(u_{(2)})$$
and
\beq
\overline{\mu}\,(u_1\ot u_2)=\mu\, (ad\,F (u_1\ot u_2))=\mu\, (ad\,
F_{(1)}(u_1)\ot ad\,F_{(2)}(u_2)).\label{ff}
\label{tw}
\eeq
Here $\mu$ is the initial product in $U(\gggg)$ and  $ad\,F^{\pm 1}$ is 
defined by 
$$ad\,X(Y)=[X,\,Y]\,\,\,{\rm  and }\,\,\,
ad\,(X_1X_2...X_p)\,(Y)=ad\,X_1\,(ad\,X_2\,(...ad\,X_p\,(Y)...).$$

The space $\ugn$ being equipped with these product, coproduct, the classical 
unit, counit and antipode becomes a twisted  Hopf algebra.
Essentially, this means that
\beq
\overline{\De}\,\overline{\mu}(u_1\ot u_2)=(\overline{\mu}\ot\overline{\mu})\, 
(\id\ot S\ot\id)(\overline{\De}(u_1)\ot \overline{\De}(u_2))
\label{twist}
\eeq
where $S=S_{\nu}=F^{-1}\sigma F$ and $\sigma$ is the flip 
($F$ and $F^{-1}$ act in the above sense by $ad\ot ad$).
Let us denote $\oH$ this twisted Hopf algebra.
We leave to the reader to verify 
that the operator $S$  satisfies 
the QYBE.

Let us observe that the $\De$-primitive elements $X\in U(\gggg)$,
i.e., such   that $\De(X)=X\ot 1+1\ot X$
(they are just the elements of the algebra $\gggg$)
are still $\overline{\De}$-primitive.
 This follows from the second relation \r{def}.

The algebra $\oH$ can be treated as the enveloping algebra of 
a generalized (or $S$-)Lie algebra defined by the deformed Lie bracket
$[\,\,,\,\,]_{\nu}=[\,\,,\,\,]\,F_{\nu}$ or in more detailed form
$$[X,Y]_{\nu}=[ad\,F_{(1)}(X), ad\,F_{(2)}(Y)].$$
An axiomatic description of such a type 
bracket is given, for example, in \cite{G4}. We will design the space 
$\gggg[[\nu]]$ equipped
with the bracket $[\,\,,\,\,]_{\nu}$ by $\gggg_{\nu}$. Its enveloping algebra 
defined naturally by 
\beq
U(\gggg_{\nu})=T(\gggg_{\nu})[[\nu]]/\{x\ot y-S(x\ot y)-[x,y]_{\nu}\}
\label{ugn}
\eeq
is filtred quadratic (more precisely, the ideal is generated by 
quadratic-linear 
elements). 

Hereafter $T(V)$ stands for the free tensor algebra of
a linear space $V$ and $\{I\}$ stands for its ideal generated by a subset 
$I\subset T(V)$.

We need also the algebra
$\ahn=U(\gggg_{\nu})_{\h}$ defined by the formula \r{ugn} but with  the bracket
$[\,\,,\,\,]_{\nu}$ replaced by $\h [\,\,,\,\,]_{\nu}$.
The algebra $\ahn$ is also filtered quadratic and moreover, possesses a 
twisted Hopf structure. Moreover, we have by construction the following

\begin{theorem} The two parameter family $\ahn$ is a flat deformation of the 
algebra
Sym$(\gggg)=\Fun(\gggg^*)$. The corresponding Poisson pencil is just \r{pp} 
where
$\{\,\,,\,\,\}_{KKS}$ is  the linear extention of
KKS bracket (Poisson-Lie one) and the  bracket $\{\,\,,\,\,\}_{R}$
is associated to the initial R-matrix.
\end{theorem}

By  passing to the quotient
$\an=\ahn/\h \ahn$ we get an S-commutative algebra which also is a flat deformation
of the algebra $\Fun(\gggg^*)$.  Let us precise that
by this we mean an algebra $\aaaa=\an$ equipped with an associative product 
$\mu:{\aaaa}^{\ot 2}\to\aaaa$ and an involutive twist
$S:{\aaaa}^{\ot 2}\to {\aaaa}^{\ot 2}$ such that $\mu\,S=\mu$ and 
$S\,\mu^{12}=\mu^{23}\,S^{12}\,S^{23}$. The latter relation signifies that 
the product
$\mu$ is S-invariant.

Now, let $\gggg$ be a simple Lie algebra. 
Then  the enveloping
algebra $\oH=U(\gggg_{\nu})$ is  isomorphic to  $U(\gggg)[[\nu]]$. 
So, we can treat it as the algebra $U(\gggg)[[\nu]]$ but equipped with a 
new coproduct
(still denoted by $\oDe$). Thus, we have equipped the algebra
$U(\gggg)[[\nu]]$ with two deformed coassociative coalgebraic structures 
converting 
it respectively into an Hopf algebra $H$ and a twisted Hopf algebra $\oH$.

However, in some sense the properties of the latter algebra are closer to
 those of the usual enveloping algebra $U(\gggg)$. 
In the first place it is due to the fact that the algebra $\oH$
possesses a generating set formed by $\oDe$-primitive elements. Moreover, for 
this 
algebra its S-commutative analogue, i.e., the algeba $\an$, is 
well defined and being equipped with the coproduct $\oDe$ is still a twisted 
Hopf algebra, as in the classical case. The passage from 
the latter algebra
to  $\ahn$ can be regarded as a twisted version of the quantization 
procedure of
the  linear Poisson-Lie bracket on $\gggg^*$ consisting in a passage 
 from the symmetric algebra
of $\gggg$ to the enveloping algebra $\ugg$. 

By means of $\oDe$-primitive elements it is not difficult to introduce
 the notion of twisted 
(or S-)vector fields: the twisted version of the Leibnitz rule for an 
involutive $S$ 
is well known. 
It is worth  noticing that the twisted vector fields are just  classical ones
 but their action on functions is deformed as follows:
\beq
\rho_{\nu}(X)\cdot a=\rho(ad\,F_{(1)}(X))\cdot\rho(F_{(2)})a, 
\,\,X\in\gggg,\,\,a\in\Fun(M) 
\label{rxa}
\eeq
where
$\rho:\gggg\to\Vect(M)$ is a representation of $\gggg$  into the space of 
vector fields on a variety $M$ extended to $U(\gggg)$.

Unfortunately, the Hopf algebra $H$  does not possess, in general, any 
generating set 
formed by $\De_F$-primitive elements (the $\De$-primitive
elements  are no longer $\De_F$-primitive). This is a reason why 
it is not so clear what is 
the natural  analogue of the
Leibnitz rule related to the quantum group $H$ (although 
some palliative forms of "quantum Leibnitz rule" can be sometimes
suggested).

Let us consider now the category of $\ugg$-modules. It can be equipped 
with the twist
$$S_{\nu}^{U,V}=(\rho_V\ot\rho_U)F^{-1}\sigma(\rho_U\ot\rho_V)F : 
U\ot V[[\nu]] \to V\ot V[[\nu]]$$  
where $\rho_U$ is the representation of $U(\gggg)$ in $U$. 
Thus, we have a twisted (symmetric monoidal in MacLane's terminilogy)
category consisting of the same objects as the initial 
one but equipped with a new transposition.

This twisted category can 
be  regarded as that of $H$-modules and that of $\oH$-modules. 
However,
the action of an element $X\in\oH$ to a tensor product of two modules $U$ and
$V$
must be defined in spirit  of the formula \r{twist} by means of the twist 
$S_{\nu}$ transposing $X_{(2)}$ and $U$ (here  $X_{(1)}\ot X_{(2)}=
\oDe(X)$). 
In particular, in this way we can deform all 
(generalized) Verma modules into twisted ones. 

Let us remark that the renormalization procedure mentioned in the introduction
(cf. also Section 3) has its twisted analogue. If in a classical case the 
map $\h\rho_{\om/\h}$
 sends $U(\gggg)$ into $\End\mom[[\h]]$ (i.e., the image does not contain
negative powers of $\h$), in a deformed case
such a property is satisfied only for an appropriated base in the deformed 
algebra.
 In the algebra
$\oH$ (which is isomorphic to $H$ as an algebra) such a base is delivered by 
$\oDe$-primitive
elements.

Let us also mention
the algebras dual to those $H$ and $\oH$. Both of them
can be treated as deformations of 
the function algebra $\Fun(G)$ on the group $G$. 
However, if the former one looks like the famous
"RTT=TTR" algebra and possesses a Hopf algebra structure, the latter one 
looks like the reflection equation (RE) algebra. 
For involutive twists it has been introduced in \cite{G1}, \cite{G4} under 
the name
of monoidal group. 
In more general setting 
RE  algebras appear as dual objects of Majid's braided groups, cf. \cite{M}. 
Majid
has also suggested a 
transmutation procedure 
converting one algebra to the other one.
RE algebras associated to twists depending on a spectral parameter were 
considered
in \cite{KS}. \\
\medskip

\noindent
{\bf 2.2} 
Let us pass  to a quasitriangular case, i.e., that 
related to the QG $\uq$ where $\gggg$ is a complex simple Lie algebra.  
In this case  there also exits 
a series $\fn$ quantizing the  R-matrix \r{rm} 
in the above sense. However, the first equation
\r{def} takes another form containing Drinfeld's associator $\Phi$ 
(cf. \cite{CP}).
 Moreover, the  corresponding twist takes the form
\beq
S=S_{\nu}=F^{-1}\sigma e^{\nu t\over 2} F \label{spc}
\eeq
where $t$ is the split Casimir.

In this case the Hopf algebra $H$ can be constructed in the same way as above.
 It is just
the famous QG $\uq$ but realized in an equivalent way as the algebra
$\ugg[[\nu]]$ equipped with the deformed coproduct $\De_F$ (we call this form
of the QG  $\uq$ its
{\em Drinfeld's realization}). However, 
the above construction of the twisted algebra $\oH$ is no longer valid because 
the product $\overline{\mu}$ defined as above is not associative 
(the associativity default is due to the Drinfeld's associator).

Nevertheless, a twisted Hopf algebra arising from  the QG $H$ exists: 
it can be obtained from $H$ by means of a transmutation procedure  which is 
 dual to that mentioned above.
In fact, this procedure does not deform the algebraic structure and transforms
the coproduct $\De_F$ into a new one $\oDe$  converting the QG into 
a "braided group".

However, this braided Hopf algebra is rather useless for us since it 
 does not apparently possess  any base of 
$\oDe$-primitive elements. 
In fact, instead of looking for an appropriated  base in $\uq$  we construct
another, complementary, algebra which possesses such a base. More precisely, we
will introduce a space $\gq$ being nothing but $\gggg$ itself equipped with an 
action $\uq\to\End \gq$ of the GQ and represent  the tensor algebra $T(\gq)$
into a q-deformed generalized Verma module $\mom^q$ with $\om=\mu\om_1$. 
Namely,
the image of the algebra $T(\gq)$ with $\mu$ expressed via $\h$ in a proper 
way 
 provides us with the quantum conterpart of the Poisson pencil \r{pp}
on the  $\cp$ type orbits (in a classic case $\h$ is proportional 
to $\mu^{-1}$ but in quantum case their relation is a little bit more 
complicated).
Hopefully, this method is valid for any symmetric orbit in $\gggg^*$ for any 
simple Lie algebra $\gggg$. 

 Note that although we do not embed the space $\gq$ into the GQ $\uq$, such 
an embedding
exists  in $sl(n)$ case  in virtue of \cite{LS}. Using this embedding the 
authors
of  \cite{LS} have introduced a version of quantum Lie $sl(n)$ bracket.

Completing this section we want to stress that in  our approach the QG 
$\uq$ play 
an auxiliary role. We use it only to describe the category which our 
quantum algebras
$\ahq$ belong to.
Let us note that in a case when such a category is related to 
nonquasiclassical twists
mentioned in footnote 2 the algebras looking like $\ahq$ can be 
constructed 
without any QG like objects. (In this case  algebras of "RTT=TTR" type can 
be introduced in the
usual way, cf. \cite{G4},  but their 
dual algebras differ drustically from the QG $\uq$.
We refer the reader to the paper \cite{AG} where an 
attempt to describe these algebras is undertaken.)

\setcounter{equation}{0}
\section{$\cp$ type orbits and their quantization by generalized Verma modules}

Let us realize now the first, classical, step of double quantization procedure
for orbits in question.

Let $\gggg$ be a simple complex Lie algebra and 
$\gggg=\hhhh\oplus\nnnn_+\oplus\nnnn_-$ be a fixed
triangular decomposition where $\hhhh$ is a Cartan subalgebra and 
$\nnnn_{\pm}$ are Borel subalgebras.
Consider a nontrivial element $\om\in\hhhh^*$ and extend it by 0 to the 
subalgebra
$\nnnn_{\pm}$. Thus, we can treat $\om$ as an element of $\gggg^*$.
Let $\coom$ be $G$-orbit of $\om$ in $\gggg^*$ where $G$ is the Lie group
corresponding to $\gggg$ acting on $\gggg^*$ by coadjoint operators and 
$$\{f,g\}_{KKS}(x)=<[df,\,dg],\,x>,\,\, x\in\coom$$
be the KKS bracket on $\coom$. 

 It is well known that the orbit $\coom$ is a closed algebraic variety in 
$\gggg^*$. 
Moreover, the space of (polynomial) functions $\aaaa=\Fun\,(\coom)$ 
can be identified with a quotient $T(\gggg)/\{I\}$ 
where $I$ is  some finit subset in $T(\gggg)$.
 
Thus, if $\coom$ is a generic semisimple orbit 
(this means that in the decomposition
$\om=\sum \mu_i \om_i$ where $\om_i$ are fundamental weights $\mu_i\not=0$ 
for any
$i$) the family $I$ consists of elements $x_ix_j-x_jx_i,\,\,1\leq i,j\leq
{\rm dim}\,\gggg$ and 
 $C_i-c_i(\om),\,\,1\leq i\leq {\rm rank}\,\gggg$ where $C_i$ are invariant
(Casimir) functions and $c_i(\om)$ are certain constants depending on $\om$.

Let us consider another example of such type orbits, namely, those 
in $\gggg^*=sl(n)^*$ of  elements  
 $\om=\mu\om_1$ or $\om=\mu\om_{n-1}$ for some $\mu\in{\bf C}$. These orbits
(or more precisely, their real compact forms in $su(n)^*$) can be identified 
with ${\bf CP}^{n-1}$. They are called $\cp$ type orbits.
It is well known that these orbits
can be described by means of a system of quadratic equations . 
An explicit form of this system follows from the structure of 
$\gggg^{\ot 2}$ as a $\gggg$-module.
Let us exhibit 
such an analysis.
\begin{proposition} Let $\gggg=sl(n),\,\,n\geq 4$.
 Then highest weights of irreducible components of $\gggg$-module
$\g^{\ot 2}$  are:
\beq
2\om_1+2\om_{n-1}, \quad \om_1+\om_{n-1}, \quad
 \om_2+2\om_{n-1}, \quad  2\om_1+\om_{n-2}, \quad
\om_2+\om_{n-2}, \quad {\rm and}\quad 0
\label{gg}
\eeq

All the irreducibles from \r{gg} occure in $\GG$ with multiplicity one
 except the irreducible with highest weight $\om_1+\om_{n-1}$ (being  highest 
weight of $\gggg$ itself)
which occures
 twice, once in the symmetric part $I_+$ of $\GG$ and once in the
skewsymmetric part $I_-$.
\end{proposition}

Note that in the $sl(2)$ case the decomposition \r{gg}
containes only the components with highest weights
$$0,\qquad 2\om_1, \qquad 4\om_1$$
 all with multiplicity one and in the $sl(3)$ case the component of 
 highest weight $\om_2+\om_{n-2}$ does not appear.

Let us denote the finite dimensional irreducible $\gggg$-module 
with highest weight $\la$ by $V_\la$. A corresponding
highest weight vector (assuming $V_\la$ to be imbedded in 
$\g^{\ot 2}$) will be denoted by $s_\la$. 
For highest weight $\om_1+\om_{n-1}$ which occures
 twice in \r{gg} we denote  $V^+_{\om_1+\om_{n-1}}$ (resp. 
 $V^-_{\om_1+\om_{n-1}}$) the component of highest weight $\om_1+\om_{n-1}$
belogning to $I_+$ (resp. $I_-$). Their  highest weight vectors will be 
designed by 
 $s^+_{\om_1+\om_{n-1}}$ (resp. $s^-_{\om_1+\om_{n-1}}$). The precise
 expressions for the corresponding highest weight vectors are presented
 in Proposition 3 with a specialization $q=1$.

Then the orbit under consideration can be defined by the following system
(here $n\geq 4$, the cases $n=2,\,\,3$ are left to the reader)
\beq
V_{\om_2+2\om_{n-1}}=0,\, V_{2\om_1+\om_{n-2}}=0,\,V^-_{\om_1+\om_{n-1}}=0
\,\,({\rm or,\,\, equivalently,}\,\, x_ix_j-x_jx_i=0\,\,\,
\forall\,\, i,j) \label{odin}
\eeq
\beq
V_{\om_2+\om_{n-2}}=0, s_0-c_0(\om)=0,\, V^+_{\om_1+\om_{n-1}}-c_1
(\om)\g=0
\label{dva}
\eeq
where $s_0=C_1$ is a generator of the trivial module 
 (Casimir element) and the constants 
$c_i(\om),\,i=0,1$ are:
\beq
c_0(\mu\om_1)=\frac{n-1}{n}{\mu}^2,\qquad 
c_1(\mu\om_1)=2\frac{n-2}{n}\mu  \label{past}
\eeq
(if we normalize $s_0,\,\,s^1_{\om_1+\om_{n-1}}$ and $s^2_{\om_1+\om_{n-1}}$ by
 \r{s5}-\r{s7} with $q=1$ and put
$s^+_{\om_1+\om_{n-1}}=s^1_{\om_1+\om_{n-1}}+s^2_{\om_1+\om_{n-1}}$). The 
last equation \r{dva} is a
symbolic form of the relation $s^+_{\om_1+\om_{n-1}}=c_1(\om) g_{1,n}$ and 
all the decendants of this
 relation.

Thus, we have ${\aaaa}=\Fun({\coom})=T(\gggg)/\{I\}$ with
the family $I\subset \bc\oplus\gggg\oplus\gggg^{\ot 2}$  
generated by the l.h.s. of the formulae \r{odin}, \r{dva}.
So, the algebra $\aaaa$ is filtred 
quadratic. (Note that this system was given in \cite{DG2} in a
nonconsistent form.)

Since the orbit $\coom$ is a symmetric space it is a spherical 
or multiplicity free
variety, i.e. in the  decomposition of the space $\Fun(\coom)$ into a 
direct sum
of irreducibles their multiplicities are at most one\footnote{Let us 
remark that
the only symmetric orbits corresponding to Lie algebra $\gggg=sl(n)$ are
$$\co_x=SL(n)/S(L(k)\times
L(n-k)),\,\,1\leq k\leq n-1$$
  whose real compact forms are Grassmanians 
(the cases $k=1$ and $k=n-1$ 
correspond to ${\bf CP}^{n-1}$).  
Symmetric orbits in $\gggg^*$ for other simple Lie algebras $\gggg$
have been classified by E.Cartan 
(cf. \cite{H}, \cite{KRR}).}. It is well known that
for the orbits of $\cp$ type
$$\Fun(\coom)\approx\bigoplus_{k=0}^{\infty}V_{k(\om_1+\om_{n-1})}.$$

Let us discuss now a way to quantize the KKS bracket well defined in the
algebras $\aaaa=\Fun(\coom)$ by means
of generalized Verma modules.

Let $K$ be the stabilizer of the point $\om\in \gggg^*$. So, $\coom=G/K$. Let
$\kkkk=$Lie$(K)$ be the Lie algebra of the group $K$ and $\pppp=\kkkk+\nnnn_+$ 
be a parabolic subalgebra of $\gggg$. Let us consider the induced 
 $\gggg$-module 
$$\mom=Ind_{\pppp}^{\gggg}\,{\bf 1}_{\om} = U(\gggg)\otimes_{U({\bf p})}
{\bf 1}_{\om}$$
 where ${\bf 1}_{\om}$ is the one
dimensional $\pppp$-module equipped with a representation 
$\rho_{\om}(x)\,e=<\om,\,x>e,\,\,x\in\pppp$ ($e$ is a generator of the module).
 The $\gggg$-module $\mom$ is usually
called {\em generalized} (in the sequel we omit this precision) 
 {\em Verma module}. 
Let us denote $\rho_{\om}:\gggg\to\End \mom$
the induced representation.

The operator algebra End$\,\mom$ is quantum object with respect
to the algebra of functions  Fun$\,(\coom)$. To give an exact meaning to this 
statement
let us introduce an associative algebra $\ah$ depending on a parameter $\h$ as 
follows.
Let us consider a map $\overline{\rho}_{\h}
=\h\rho_{\om/\h}:\gggg\to \End \mom[[\h]]$,
extend it naturally to $T(\gggg)$
and introduce the algebra
$\ah$ as subalgebra of End$\,\mom[[\h]]$ being, by definition, the image
 $\overline{\rho}_{\h}(T(\gggg))$. 

\begin{proposition}
The algebra $\ah$ is a flat deformation of that ${\aaaa}=Fun(\coom)$  
and the corresponding Poisson bracket  is just the KKS one. 
\end{proposition}

 This statement is valid for any simple Lie algebra and for any ss orbit.
We  demonstrate it for the orbits of ${\bf CP}^n$ type where all the 
 calculations can be easily done. 
In fact we will see that in this case the algebras $\aaaa$ and  $\ah/\h\ah$ are
isomophic as $\gggg$-modules (they are consisting of the same irreducibles
with multiplicity one).

Let us first consider the  
finite dimensional $sl(n)$-modules $V_{\om},\,\,\om=\mu\om_1,\,\,
\mu\in{\bf Z}_+$
where ${\bf Z}_+$ stands for the set of nonnegative integers.
Such a module can be naturally identified with symmetric power of the 
vector fundamental space $V_{\om_1}$ (in the $sl(2)$ case the factor 
$\mu$ is just spin
of the module).
Its dimension is equal to 
$\left(\begin{array}{c} \mu+n-1\\n-1\end{array}\right)$.

Let us fix in the space $V_{\om}$
the base 
\beq
\mid m_1,\ldots ,m_n> =x_1^{m_1}\cdots x_n^{m_n},\,\, 
\sum m_i = \mu. \label{base}
\eeq

Let $h_i\in\hhhh,\,e_i\in\nnnn_+,\,f_i\in\nnnn_-,\,\,1\leq i\leq n-1$
be a standart Chevelley base in the Lie algebra $sl(n)$. The elements
$h_i,\,e_i,\,f_i$ act in the 
module $V_{\om}$ as first order differential operators 
\beq
e_i=x_i\frac{\partial}{\partial x_{i+1}}\ ,\qquad
f_i=x_{i+1}\frac{\partial}{\partial x_{i}}\ ,\qquad
h_i=x_i\frac{\partial}{\partial x_{i}}-
x_{i+1}\frac{\partial}{\partial x_{i+1}}\ .
\label{3a}
\eeq

In the base \r{base} the operators 
\r{3a} look like 
\bea
e_i\mid m_1,\ldots ,m_n>&=&m_{i+1}\mid m_1,\ldots,m_i+1,m_{i+1}-1,\ldots,m_n>
 \ ,\nn\\
f_i\mid m_1,\ldots ,m_n>&=&m_{i}\mid m_1,\ldots,m_i-1,m_{i+1}+1,\ldots,m_n>
 \ ,\label{3b}\\
h_i\mid m_1,\ldots ,m_n>&=&\left(m_i-m_{i+1}\right)\mid
m_1,\ldots,m_i,m_{i+1},\ldots,m_n>\ . \nn
\eea

It is well known that $sl(n)$-module $\End V_{\mu\om_1}$, $\mu\in{\bf Z}_+$
is isomorphic to the
following
 multiplicity free direct sum:
\beq
\End V_{\mu\om_1}\approx\bigoplus_{k=0}^{\mu} V_{k(\om_1+\om_{n-1})}.
\label{End}
\eeq

Let us pass now to the  Verma module 
$\mom,\,\om=\mu\om_1,\,\mu\in
{\bf C}$. Similary to the above finite dimensional modules $\vom$
 it possesses the following  base 
\beq
\mid m_1,\ldots, m_n>,\qquad \sum m_k = \mu, 
\quad m_k\in {\bf Z}_+, \,\,\,k= 2,\ldots, n.
\label{base1}
\eeq
Thus, the elements of the base \r{base1} are labeled
by the vectors $(m_2,\ldots,m_n),\,\,m_k\in {\bf Z}_+$ .
 The action of $sl(n)$ on $M_{\mu\om_1}$ is
 given by the formulae \r{3a}-\r{3b} as well.

The formula \r{End} must be modified as follows 
\beq
\Im\rho_{\om}(T(\gggg))\approx\bigoplus_{k=0}^{\infty} 
V_{k(\om_1+\om_{n-1})}.
\label{fff}
\eeq
We have replaced $\End V_{\mu\om_1}$ by 
$\Im\rho_{\om}(T(\gggg))$ since for infinite dimensional modules the map 
$\rho_{\om}$ is no longer surjective.

Let us now go back  to the algebra 
$\ah=\Im\overline{\rho}_{\h}(T(\gggg))\subset 
\End \mom[[\h]]$. By using the decomposition \r{fff} it is easy to show that
$\ah$ is a flat deformation of the algebra $\aaaa$. 
The map $\overline{\rho}_{\h}$ sends the Chevalley generators 
into operators acting with respect to the formulae 
 \r{3b}  but with $m_2,...,m_n$ 
replaced by $\h m_2,...,\h m_n$. The commutators
between the images of the Chevalley generators
 are  those in $sl(n)$ multiplied by $\h$. This implies that the 
corresponding Poisson bracket is equal to the KKS one.

Let us represent now the algebra $\ah$ as a quotient 
$$\ah=\Im\overline{\rho}_{\h}(T(\gggg))=
T(\gggg)[[\h]]/\Ker\overline{\rho}_{\h}.$$
In the case under consideration $(\gggg=sl(n),\,\om=\mu\om_1)$ this 
quotient is also a quadratic algebra.
More precisely, the ideal $\Ker\overline{\rho}_{\h}$ is generated 
by a finite family
$I_{\h}\in{\bf C}\oplus\g\oplus\g^{\ot 2}$, looking like that
defined by the l.h.s. of  \r{odin}-\r{dva} but  with some
evident modifications: the elements $x_ix_j-x_jx_i$
must be replaced by those
$x_ix_j-x_jx_i-\h[x_i,x_j]$ and the factors $c_i(\om),\,\,i=0,1$
 must be deformed to
those $c_i(\om,\h)$ depending on $\h$ (with $c(\om,\h)=c(\om)$ mod $\h$).
Namely, with $s_0$ and $s^+_{\om_1+\om_{n-1}}$ normalized as in \r{past} 
we have
$$c_0(\mu\om_1,\h)=\frac{n-1}{n}\mu(\mu+n\h),\qquad c_1(\mu\om_1,\h)=
\frac{n-2}{n}(2\mu+n\h)\ .$$

\setcounter{equation}{0}
\section{Braided algebras}

Our next aim is to braid the above quantization procedure. Let us begin with a 
description of the space $\gq$ mentioned in Section 2.

Let $\Uq$ be the quantum enveloping algebra corresponding to
 $sl(n)$. In Chevalley  generators $e_i, f_i, h_i$, $ i=1, \ldots ,n-1$
 it could be described by the relations 
\beq
[h_i,e_i]=2e_i, \qquad  [h_i,f_i]=-2f_i,
\label{0}
\eeq
\beq
[h_i,e_{i\pm 1}]=-e_{i\pm 1},\qquad [h_i,f_{i\pm 1}]=f_{i\pm 1}
\eeq
\beq
[h_i,e_j]=[h_i,f_j]=0,\quad \mid i-j\mid >1,
\eeq
\beq
[e_i,f_j]=\delta_{i,j}\frac{q^{h_i}-q^{-h_i}}{q-q^{-1}},
\label{1}
\eeq
\beq
e_i^2e_{i\pm 1}-[2]_qe_ie_{i\pm 1}e_i+e_{i\pm 1}e_i^2=0,
\label{2}
\eeq
\beq
f_i^2f_{i\pm 1}-[2]_qf_if_{i\pm 1}f_i+f_{i\pm 1}f_i^2=0,
\label{3}
\eeq
with 
$$ [n]_q=\frac{q^n-q^{-n}}{q-q^{-1}} \quad{\rm and}\quad q^{\alpha h_i}=
\exp(\nu \alpha h_i)\ .$$

We choose a comultiplication map as follows:
\beq
\Delta h_i =h_i\ot 1+1\ot h_i,\quad \Delta e_i =e_i\ot 1+q^{-h_i}\ot e_i,
 \quad \Delta f_i =1\ot f_i+f_i\ot q^{h_i}\ .
\label{4}
\eeq
 Then the antipode has a form:
$$s(h_i)=-h_i,\quad s(e_i)=-q^{h_i}e_i,\quad s(f_i)=-f_iq^{-h_i}\ .$$

Let $\gq$ be a $q$-analogue of an adjoint representation of Lie algebra $sl(n)$
 on itself, i.e., $\gq$ is $(n^2-1)\ $-dimensional $\Uq$--module with
 highest weight $\om_1 +\om_{n-1}$. We want to describe the 
 action of the QG $\Uq$ to $\gq$ explicitely 
in a fixed base of $\gq$. We denote
further
 this action by $\ad=\ad_q$. Namely, the vector space $\gq$ is generated by the
elements $g_{i,j}$, $i,j = 1,\ldots ,n,\  i\not = j$ and 
$t_i, i = 1,\ldots ,n-1\ .$
 An action of Cartan elements coincides with the classical one:
$$\ad h_i(t_k) =0\ ,$$ 
\beq
\ad h_i (\gkl)=\left(\delta_{i,k}-\delta_{i,l}-\delta_{i+1,k}+\delta_{i+1,l}
\right) \gkl
\label{5}
\eeq
Nontrivial matrix coefficients of the action of Chevalley genarators 
$e_i$ and $f_i$ of $\Uq$  look as follows:
\beq
 \ad e_i(g_{a,i})= -g_{a,i+1}\ , \qquad 
 \ad e_i(g_{i+1,a})= g_{i,a}\ ,\qquad a\not =i,i+1\ ,
\label{6}
\eeq
 $$ \ad e_i(g_{i+1,i})= t_{i}\ , \qquad 
 \ad e_i(t_{i})= -[2]_qg_{i,i+1}\ ,\qquad
\ad e_i(t_{i\pm 1})= g_{i,i+1}\ ,$$
 $$ \ad f_i(g_{a,i+1})= -g_{a,i}\ , \qquad 
 \ad f_i(g_{i,a})= g_{i+1,a}\ ,\qquad a\not =i,i+1\ ,$$
$$
\ad f_i(g_{i,i+1})= -t_{i}\ , \qquad 
 \ad f_i(t_{i})= [2]_qg_{i+1,i}\ ,$$
\beq
\ad f_i(t_{i\pm 1})= -g_{i+1,i}\ .
\label{7}
\eeq

So, we get the matrix coefficients of this action from
 the classical ones replacing the coeffitient $2$ by its 
$q$-analogue $[2]_q = q+q^{-1}$.

Since $\gq$ is a $\Uq$--module, the tensor algebra $\Tgq$ can be equipped with 
a $\Uq$--invariant product in sense of \r{1.1}. In the what 
follows the algebra
$T(\gq)$ and all its $\uq$-invariant quotients will be called {\em 
braided} ones.

In fact, the braided algebra $T(\gq)$ is "too big" for us. We are rather
 interested  in its quotient over the kernel of map sending this algebra
into $\End\,\mom^q$ where $\mom^q$ is a q-analogue of the above Verma modules
with $\om=\mu\om_1$. Namely this quotient  with $\mu$ properly expressed 
 via the parameter $\h$ plays the role of our "double quantum"
object $\ahq$. 

Let us describe the  mentioned kernel. To do this
 we need a decomposition of the $\Uq$-module $\gqq$ into a direct
sum of irreducibles.

\begin{proposition}
The formulas below describe all highest weight vectors of reducibles 
in $\Uq$-module
 $\gqq\,\, (n\geq 4)$:
\bea
s_{2\om_1+2\om_{n-1}}&=&g_{1,n}\ot g_{1,n}\ ,
\label{s1}\\
s_{2\om_1+\om_{n-2}}&=&g_{1,n}\ot g_{1,n-1}-q^{-1}g_{1,n-1}\ot g_{1,n}\ ,
\label{s2}\\
s_{\om_2+2\om_{n-1}}&=&g_{1,n}\ot g_{2,n}-q^{-1}g_{2,n}\ot g_{1,n}\ ,
\label{s3}\\
s_{\om_2+\om_{n-2}}&=&qg_{1,n}\ot g_{2,n-1}+q^{-1}g_{2,n-1}\ot g_{1,n}
-g_{1,n-1}\ot g_{2,n}-g_{2,n}\ot g_{1,n-1}
\ ,
\label{s4}\\
s^1_{\om_1+\om_{n-1}}&=&
g_{1,2}\ot g_{2,n}+q g_{1,3}\ot g_{3,n}+\ldots +
q^{n-3} g_{1,n-1}\ot g_{n-1,n}+\nn\\
&+&q^{-2}\sum_{k=1}^{n-1}\frac{[n-k]_q}{[n]_q}t_k\ot g_{1,n}-
q^{n-2}\sum_{k=1}^{n-1}g_{1,n}\ot \frac{[k]_q}{[n]_q}t_k \ ,
\label{s5}\\
s^2_{\om_1+\om_{n-1}}&=&
g_{2,n}\ot g_{1,2}+q^{-1}g_{3,n}\ot g_{1,3}+\ldots +
q^{-n+3}g_{n-1,n}\ot g_{1,n-1}-\nn\\
&-&q^{1-n}\sum_{k=1}^{n-1}\frac{[k]_q}{[n]_q}t_k\ot g_{1,n}+
q\sum_{k=1}^{n-1}g_{1,n}\ot \frac{[n-k]_q}{[n]_q}t_k \ ,
\label{s6}\\
s_0&=&\sum_{i,j=1,i\leq  j}^{n-1}\frac{[i]_q[n-j]_q}{[n]_q}t_i\ot t_j +
\sum_{i,j=1,i>  j}^{n-1}\frac{[j]_q[n-i]_q}{[n]_q}t_i\ot t_j\nn\\
&+& q\sum_{i<j}^{}q^{j-i}\gij\ot g_{j,i}+
q^{-1}\sum_{i>j}^{}q^{j-i}\gij\ot g_{j,i} \ .
\label{s7}
\eea
\end{proposition}

Because of multiplicity in this decomposition of the  highest weight
$\om_1+\om_{n-1}$ component it is not so clear 
what are natural q-analogues $\Ipmq$ of symmetric
$I_+$ and skewsymmetric $I_-$ components in $\gqq$ (exept the $sl(2)$ case)
\footnote{For
other simple Lie algebras $\gggg$ such natural q-analogues $\Ipmq$ exist 
since $\gggg^{\ot 2}$ is 
multiplicity free as a $\gggg$-module but deformations $T(\gggg)/\{I_{\pm}\}\to
T(\gq)/\{\Ipmq\}$ are not flat. In $sl(n)$ case one can split $\gqq$ into a 
direct sum
of components $\Ipmq$ in such a way that these deformations are flat.
It is shown in \cite{Do} by means of an embedding $\gq\to\uq$ which is 
slightly different from that considered in \cite{LS}.}.

Let us consider a way to introduce  a decomposition $\gqq=I^q_+\oplus I^q_-$ 
arising from an 
 operator $\tS$ discussed in \cite{DS} and \cite{DG2}.
In Drinfeld's realization of the QC $\uq$ the operator $\tS$ 
is defined by the formula \r{spc} but without the factor 
$e^{\nu t\over 2}$. So, it is evident that this operator is involutive. 
Moreover, 
being restricted on $\gqq$
it has the same
eigenspaces as the YB operator $S$ has but with  eigenvalues $\pm 1$.
Namely, to pass from $S$ to $\tS$ we must replace
the eigenvalues of $S$ close to 1 (resp. -1) by 1 (resp -1) assuming that 
$\vert
q-1\vert\ll 1$.

We complete this Section with describing the action of $\tS$ on the 
isotypical component
of highest weight $\om_1+\om_{n-1}$
(its action on other components of $\gqq$  containes no new information 
for us). We use this computation in
the last Section.
 
To do this we need a partial information on the quantum universal
R-matrix in an $sl(n)$ case.

It is well known that the universal $R$-matrix $\R$ for the
 algebra $U_q(\g)$ can be presented by
\bn
\R= \R_0\cdot q^{\sum c_{i,j}h_i\otimes h_j}\ ,
\label{rr1}
\ed
 where $(c_{i,j})$ is the matrix  inverse to the Cartan matrix of $\g$ and
$\R_0$ belongs to a tensor product of quantized enveloping algebras
 of nilpotent subalgebras $\n_\pm$ of $\g$:
$\R_0\in U_q(\n_+)\otimes U_q(\n_-)$. Moreover, $\R_0=1$ mod 
$\n_+U_q(\n_+)\otimes U_q(\n_-)$.

In  $sl(n)$ case formula (\ref{rr1}) has especially
 simple form after embedding of $\R$ into $U_q(gl_n)\otimes U_q(gl_n)$:
\bn
\R= \R_0\cdot q^{\sum_{i=1}^n \varepsilon_i\otimes \varepsilon_i -
\frac{1}{n}(\sum_{i=1}^n \varepsilon_i)\otimes (\sum_{i=1}^n \varepsilon_i)}
\ ,
\label{rr2}
\ed
 where $h_i=\va_i-\va_{i+1}$, $i=1,\ldots ,n-1$.

Let us compute the expressions
$S(\so^i)$, $i=1,2$, where $S=\sigma (ad \otimes ad) \R$
 is the image of the universal $R$-matrix in tensor square of adjoint
 representation multiplied by the flip. The commutativity of
 $S$ and $\Delta(x)$ imlplies that
$S(\so^i)=\sum_j a_{i,j}\so^j,\quad i,j=1,2$
 for some constants $a_{i,j}$.

The space $\gq$ can be decomposed into three parts:
$$\gq=\g_++{\bf t}+\g_-$$
where $\g_+$ is generated by the vectors $g_{i,j}, i<j;\ $
 $\g_-$ is generated by the vectors $g_{i,j}, i>j;\ $ and
 ${\bf t}$ is generated by elements $t_i$. Their crucial properties are:
 $$U_q(\n_\pm){\bf t}\subset \g_\pm \,\,{\rm and} \,\, 
U_q(\n_\pm)\g_\pm\subset \g_\pm.$$
One can observe from the explicit expressions for $\so^{i}$ that
\bn \so^1=q^{-2}\sum_{k=1}^{n-1}\frac{[n-k]_q}{[n]_q}t_k\otimes g_{1,n}
+\sum x_i\otimes y_i \ ,
\qquad x_i\in\g_+\ ,
\label{rr3}
\ed
\bn
\so^2=q\sum_{k=1}^{n-1}\frac{[n-k]_q}{[n]_q}g_{1,n}\otimes t_k
+\sum u_i\otimes v_i \ ,
\qquad v_i\in \g_+\ .
\label{rr4}
\ed

Due to (\ref{rr2}) we have
$$S(\so^1)=q^{-2}\sum_{k=1}^{n-1}\frac{[n-k]_q}{[n]_q}t_k\otimes g_{1,n}
+\sum {x'}_i\otimes {y'}_i \ ,
\qquad {x'}_i\in\g_+\ ,,
$$
 which means, in virtue of (\ref{rr4}), that
\bn S(\so^1)=q^{-3}\so^2\ .
\label{rr5}
\ed
Analogously,
\bn
S(\so^2)=q^{3-2n}\so^1\ .
\label{rr6}
\ed

Formulas (\ref{rr5}) and (\ref{rr6}) show that the operator $S$
 is diagonal on isotypical component $V^q_{\omega_1+\omega_{n-1}}\oplus 
V^q_{\omega_1+\omega_{n-1}}$
 (as well as in the whole space $\gqq$), has on it eigenvalues
$\pm q^{-n}$, and the corresponding eigenvectors are 
\beq
s_\pm= q^{2-n}\so^1\pm q^{-1}\so^2.
\label{sm}
\eeq
This result is true for $n\geq 3$.
The case $n=2$ is left for the reader (here $s_+=0$).

\setcounter{equation}{0}
\section{Braided modules}

\begin{definition} We say that  $M$ is a braided $\Tgq$--module (or, simply
braided module), if
$M$ is equiped with a
 structure of $\Uq$--module  and of $\Tgq$--module, and these structures are 
 related as 
\beq
u\cdot (gm)=\left( u^{(1)}\cdot g\right)
\left( u^{(2)}\cdot m\right)
\label{9}
\eeq
for any $u \in \Uq$, $g \in \Tgq$ and $ m\in M$.
\end{definition}

The braided algebra $\Tgq$ together with the  category 
  of its braided representations can be described also in a language
 of intertwining operators. Let $M$ be a $\Uq $-module. Then, by definition,
the (second type) intertwining operator $\Psi^{\g_q}$ is a $\Uq$-morphism.
\beq
\Psi^{\g_q}\ : \g_q\ot M\rightarrow M
\label{2a1}
\eeq
 The components $\Psi^{\g_q}_a\ :M \rightarrow M$ are definined via 
 fixing a base $g_a$ in $\gq$:
\beq
\Psi^{\g_q}_a(m)= \Psi^{\g_q}(g_a\ot m)\ .
\label{2a2}
\eeq
Thus, $M$ is a braided $\Tgq$-module if and
only if there exists an action of intertwining operator $\Psi^{\g_q}$ on $M$ 
in the above sense.

Our next aim is to perform an explicit construction of certain braided modules.
More precisely, we will define a $\Uq$-morphism
\beq
 T(\gq)\to \End \vom^q,\,\om=\mu\om_1,\,\,\mu\in{\bf Z}_+ \label{rrr}
\eeq
where $\vom^q$ is a q-deformed finite dimensional
module. After that we will extend this construction  to the q-deformed Verma
modules. 

Since the module $\End\,\vom$ is multiplicity free the component isomorphic 
to $\gggg$ 
 is represented in it only once. The same is true for the $\Uq$-module
$\vom^q$. This enables us to define the map \r{rrr} in a unique way up to a 
factor assuming it to be a
$\Uq$-morphism. Moreover, the space $\End\,\vom$ does not contain any 
component with highest
weights $\om_2+2\om_{n-1}$ and $2\om_1+\om_{n-2}$ (it is also true
in a q-deformed case). 

The $\uq$-modules possessing these two properties were called in \cite{G5}
braided, here we use this term in more general sense.

Let us describe now the map \r{rrr} explicitly. 
Since the algebra $\Tgq$ is generated by the space $\gq$ and $\Uq$ 
is generated 
 by Chevalley base, it suffices to  ensure the relation \r{9} for 
 $u =h_i, e_i, f_i$ and for $g\in \gq$. Below we write down these relations
 using the following traditional notation. Let  $M$ be a $\Uq$--module and for
  $x \in M$ be an eigenvector of
the action of Cartan subalgebra $\hhhh$ of $\Uq$. Then we denote its eigenvalue
by 
 $\la(x) \in \hhhh^*$ such that 
 $$
h_i(x)=\left(\ve_i-\ve_{i+1}, \la(x)\right) \cdot x\ ,\qquad
 \ve_i={\rm diag}\left( 0,\ldots ,0,
 1,¯,\ldots ,0\right)\ ,  {\mbox{1 at {\it i}-th place}} . 
$$
 For instance, the weights $\la(g)$ 
of the representation ad in a space $\gggg_q$
 coincide  with the classical ones:
$$ \la(\gij)=\ve_i-\ve_j,\qquad \la(t_i)=0\ .$$
and the action of the Cartan elements in $\gggg_q$ is given by the relation
$\ad{h_i}(g)= \left(\ve_i-\ve_{i+1}, \la(g)\right)\cdot g$.
\begin{proposition} Let $M$ be a finite dimensional $\Uq$--module. Then $M$ is
braided $\Tgq$--module if and only if for any $g\in \gq$  the following 
relations
hold:
\beq
[h_i,g]=\ad{h_i}(g)\,,
\label{10a}
\eeq
\beq
 [e_i,g]_{q^{-\left(\ve_i-\ve_{i+1},\la(g)\right)}}=
\ad{e_i}(g)\ ,
\label{10b}
\eeq
\beq
[f_i,g]=\ad{f_i}(g)\cdot q^{h_i}\ .
\label{10}
\eeq
Here 
$$[a,b]_q = ab-qba$$ 
and all brackets are understood in operator sense.
\end{proposition}

{\bf Proof} It suffices to say that 
for any finite dimensional $\Uq$-module $M$  we identify 
$$\End M=M \ot M^*$$
 as left $\Uq$-modules where an action of $\Uq$ to $M^*$ is defined 
 by means of an antipode $s$: 
$$( v, u\cdot\xi)=(s(u)\cdot v,\xi),\qquad v\in M,\; \xi\in M^*\;
 u\in\Uq\ .$$
 The rest is a substitution of \r{4} in \r{9}.

Let $V_{\om_1}^q$ be the first (vector) fundamental representation of the
 algebra $\Uq$.
Similary to classical case let us consider the  irreducible 
finite dimensional representations 
$\Vm^q, \mu \in{\bf Z}_+$ of $\Uq$ with  highest weights $\mu\om_1$
 (in what follows we will omit $q$).

One can easily check that the operators $e_i, f_i, h_i \in{\rm End\ }\Vm$, 
 whose nontrivial matrix elements are desribed in \r{2.3}  satisfy the 
 relations \r{0}-\r{3} and thus define an action of the algebra $\Uq$ in vector
 space $\Vm$:
\bea
e_i\mid m_1,\ldots ,m_n>&=&
[m_{i+1}]_q\mid m_1,\ldots,m_i+1,m_{i+1}-1,\ldots,m_n>
 \ ,\nn\\
f_i\mid m_1,\ldots ,m_n>&=&
[m_{i}]_q\mid m_1,\ldots,m_i-1,m_{i+1}+1,\ldots,m_n>
 \ ,\label{2.3}\\
h_i\mid m_1,\ldots ,m_n>&=&\left(m_i-m_{i+1}\right)\mid
m_1,\ldots,m_i,m_{i+1},\ldots,m_n>\ . \nn
\eea

\begin{proposition} There is a unique structure (up to a multiplicative
 constant $\alpha \in {\bf C}$) of braided $\Tgq$--module
 on $\Uq$--module $\Vm$. 

The action of generators of $\Tgq$ in braided module
$\Vm$ is:
$$
\gij \mid  m_1,\ldots ,m_n>
 =\a(\mu)\  
q^{ j+(m_1+\ldots +m_i)-(m_j+\ldots +m_n)}\cdot$$
\beq
\cdot 
[m_{j}]_q\mid m_1,\ldots,m_i+1,\ldots ,m_{j}-1,\ldots,m_n>
 \label{2.4}
\eeq
$$
{\rm for\ } i<j
$$
 $$
g_{j,i}\mid  m_1,\ldots ,m_n>=\a(\mu)\ 
q^{ i-1+(m_1+\ldots +m_{i-1})-(m_{j+1}+\ldots +m_n)}\cdot
$$
\beq
\cdot 
[m_{i}]_q\mid m_1,\ldots,m_i-1,\ldots , m_{j}+1,\ldots,m_n>
\label{2.5}
\eeq
$$
{\rm for\ } i<j
$$
 $$
t_{i}\mid  m_1,\ldots ,m_n>=\a(\mu)\ 
q^{ i+(m_1+\ldots +m_{i-1})-(m_{i+2}+\ldots +m_n)}\cdot$$
\beq
\cdot
\frac{\left( [2]_qq^{m_i-m_{i+1}}-q^{m_i+m_{i+1}+1}-
q^{-m_i-m_{i+1}-1}\right)}{q-q^{-1}} 
\mid
m_1,\ldots,m_n>\ .
\label{2.6}
\eeq
\end{proposition}

{\bf Proof} The proof goes by induction on the rank $n$.
Let us start from the $U_q(sl(2))$ case. In this case the relation \r{10}
implies that 
\beq
[g_{2,1}, f_1]=0\ .
\label{g21}
\eeq      
From
the commutation relation \r{10a} of $g_{2,1}$ with Cartan element we see that
 in addition $g_{2,1}$ has the same matrix structure as $f_1$ and thus 
these two
 operators are proportional
to each other. Applying twice the relation \r{10b} we get the 
 description of the operators $t_1$ and $g_{1,2}$. Finally we check that
all the relations \r{10a}--\r{10} are satisfied. 

The passage to $U_q(sl(3))$
 looks as follows. We know from $sl(2)$ case that the operator
 $g_{2,1}$ has a form 
\beq
 g_{2,1}=\a(m_3)f_1\ .
\label{am3}
\eeq
We wish to find out the normalization constant $\a(m_3)$. Applying the 
following particular cases of \r{10b} and of \r{10}:
$$g_{3,1}=[f_2,g_{2,1}]q^{-h_2}, \qquad 
 g_{3,2}=-[e_1,g_{3,1}]$$
 to the ansatz \r{am3}, we get the description of the operator $g_{3,2}$,
 depending on a choise of $\a(m_3)$. But we know again that 
$$
g_{3,2}=\a(m_1)f_2\ .
$$
This gives a reccurence equation on $\a(m_3)$ which unique (up to a constant
 factor) solution is $\a(m_3) = q^{-m_3}$. Then we get from \r{10} the
description of others generators of $\gq$ in $sl(3)$ case. The general 
induction 
step is similar.

Let us pass now to q-deformed Verma modules. Let  $\mom=\mom^q$ be such a 
module.
It also possesses a base 
labeled by $(m_2,\ldots,m_n),\,m_i\in {\bf Z}_+$.
For any $\mu\in {\bf C}$
there exists a $\Uq$-invariant map 
\beq
\orom:T(\gq)\to \End \mom,\,\om=\mu\om_1. \label{oro}
\eeq
This map is also defined (uniquely up to a factor)
 by the formulae \r{2.4}--\r{2.6}.

We are interested now in the ideal $\Ker\orom$. It is also generated by its
quadratic part  $\iqma=\Ker\orom\bigcap(\bc\oplus
\gq\oplus\gqq)$. To describe this quadratic part we consider the images 
of the  
highest weight elements $s_0$, $s^1_{\om_1+\om_{n-1}}$ and 
$s^2_{\om_1+\om_{n-1}}$ in $\gqq$ with respect to the map $\orom$.

 The operator $s_0$ is a scalar (we omit the symbol $\orom$) :
\beq
 s_0=\a(\mu)^2q^n\frac{[n-1]_q}{[n]_q}[\mu]_q[\mu+n]_q\ {\rm Id}\ ,
\label{laplace}
\eeq
and the operators $s^1_{\om_1+\om_{n-1}}$ and $s^2_{\om_1+\om_{n-1}}$ 
 are proportional to the operator $g_{1,n}$:
\bea
s^1_{\om_1+\om_{n-1}}&=&\a(\mu)q^{n-2}\frac{[n-1]_q
[\mu+n]_q -[\mu]_q}{[n]_q}\ g_{1,n}\ ,
\label{g1}\\
s^2_{\om_1+\om_{n-1}}&=&\a(\mu)q\frac{[n-1]_q
[\mu]_q -[\mu+n]_q}{[n]_q}\ g_{1,n}\ .
\label{g2}
\eea

Finally, we have the following description of $\iqma$.

\begin{proposition}  The subspace $\iqma \subset \left(
{\bf C}\oplus\gq\oplus\gqq\right)$ 
is a $\Uq$-module generated by

(i)  the highest weight vectors $s_{2\om_1+\om_{n-2}}$,
$s_{\om_2+2\om_{n-1}}$, $s_{\om_2+\om_{n-2}}$

(ii)  the following combinations of highest weight vectors:
\bea
s^1_{\om_1+\om_{n-1}}&-&\a(\mu)q^{n-2}\frac{[n-1]_q
[\mu+n]_q -[\mu]_q}{[n]_q}\ g_{1,n}\ 
\label{gg1}\\
s^2_{\om_1+\om_{n-1}}&-&\a(\mu)q\frac{[n-1]_q
[\mu]_q -[\mu+n]_q}{[n]_q}\ g_{1,n}\ .
\label{gg2}\\
 s_0&-&\a(\mu)^2q^n\frac{[n-1]_q}{[n]_q}[\mu]_q[\mu+n]_q\ \cdot\ {1}\ ,
\label{gg3}
\eea
\end{proposition}

Therefore for the elements $s_\pm$ defined by \r{sm} we have the 
following formula
\bea
s_\pm&=&\a(\mu)\frac{[n-1]_q\mp 1}{[n]_q}([\mu+n]_q\pm [\mu]_q)\ g_{1,n}\ . 
\label{sm1}
\eea

\setcounter{equation}{0}
\section{The algebra $\ahq$ and quantum $\cp$ type orbits}

Let us define now the two parameter algebra $\ahq$ using the 
results of the previous
Sections. To do this we must express $\mu$ via $\h$ and choose 
the factor $\a(\mu)$ in a proper way.
In the classical case $(q=1)$ by setting $\h=\mu^{-1},\,\,
\a(\mu)=\a_0\h$ we get an algebra  which differs from the above 
algebra $\ah$ by a renormalization of the parameter. 
Thus, the algebra $\ah/\h\ah$ is just 
function algebra on the corresponding orbit (labeled by $\a_0$).

In the quantum case $(q\not=1)$ we suppose that $\vert q\vert\not=1$. 
This 
condition is motivated by our desire to have $[\mu]_q\to\infty$ as
$\mu\to\infty$. 

Let us set 
\beq
\alpha(\mu)=\frac{\alz}{ [\mu]_q}\,\,\,{\rm and}\,\, \,
\frac{[\mu+n]_q}{[\mu]_q}=\gamma(q)+\h,\label{mm}
\eeq
where $\gamma(q)=q^n$ if $\vert q\vert>1$ and $\gamma(q)=q^{-n}$ if 
$\vert q\vert<1$. We use \r{mm} as the definition of the parameter $\h$.

Then the  elements  \r{gg1}-\r{gg3} become
\bea
s^1_{\om_1+\om_{n-1}}&-&\alz q^{n-2}\frac{[n-1]_q
(\gamma(q)+\h) -1}{[n]_q}\ g_{1,n}\ 
\label{gg11}\\
s^2_{\om_1+\om_{n-1}}&-&\alz q\frac{[n-1]_q
 -(\gamma(q)+\h)}{[n]_q}\ g_{1,n}\ .
\label{gg21}\\
 s_0&-&\alz^2q^n\frac{[n-1]_q}{[n]_q} (\gamma(q)+\h)\ \cdot\ {1}\ .
\label{gg31}
\eea

Meanwhile, the element $s_-$ defined by the formula \r{sm1}  takes the form
\bea
s_-&=&\alz\frac{[n-1]_q+1}{[n]_q}((\gamma(q)+\h)-1)\ g_{1,n}\,. \label{sm2}
\eea

Let us introduce now the algebra $\ahq$ as the quotient of $T(\gq)$
by the ideal generated by the elements listed in 
 Proposition 6 (i), the elements \r{gg11}-\r{gg31} and all their
descendants.  Expressing $\mu$ via $\h$ and sustituting it in the formulae 
\r{2.4}-\r{2.6}
we can realize the algebra $\ahq$ as some subalgebra in $\End \vom^q[[\h]]$ 
with a fixed $\om$ 
($\alpha(\mu)$ is assumed to be expressed via the formula \r{mm}).

This operator realization of the algebra $\ahq$ implies that deformation 
${\aaaa}\to{\ahq}$
 is flat. In fact, it suffices to note that the algebra $\ahq$ contains all 
compoments 
$V_{k(\om_1+\om_{n-1})},\,\,k=0,1,2,...$ 
This follows from the fact that the images of the elements
$g_{1,n}^k$ are not trivial operators for any $k=0,1,2,...$ (recall that
$\vert q \vert\not=1$ and therefore $q$ is not any root of the unity).

The arguments analogous to  Nakayama lemma (cf. \cite{AM}) 
show that the algebra $\aq=\ahq/\h\ahq$ also contains all components
$V_{k(\om_1+\om_{n-1})}$ and therefore the deformation ${\aaaa}\to{\aq}$
 is flat. 

Our next aim is to verify that the Poisson pencil corresponding to the algebra
$\ahq$ is just that \r{pp} with R-matrix 
\r{rm}. To do this it suffices to compute the brackets corresponding to one 
parameter
deformations ${\aaaa}\to {\aq}$ and ${\aaaa}\to {\ahq/(q-1)\ahq}$. It is 
easy to 
see that the algebra $\ahq/(q-1)\ahq$ is just that discussed in the
begining of this
Section. Therefore the corresponding Poisson bracket is proportional to the 
KKS one.

Consider now the algebra $\aq$. 
Let us remark that this algebra differs from 
analogous one parameter algebras from \cite{DS} and \cite{DG2}.
The latter algebras were $\tS$-commutative where the operator $\tS$ is 
defined in Section 4
and the  algebra
$\aq$ is no longer $\tS$-commutative.
Instead of the realation $s_-=0$ taking place in an $\tS$-commutative
algebra we have now \r{sm2} with $\h=0$.
In \cite{DG2} it has been shown that Poisson bracket corresponding
to $\tS$-commutative algebra on a symmetric orbit is proportional to
the R-matrix one.

The $\tS$-commutativity default of the  algebra $\aq$ is mesured by the 
r.h.s. of
the formula \r{sm2}. In the quasiclassical limit this term gives rise
to a contribution proportional to the KKS bracket. This completes the proof.

In this connection the following question arises: what algebra of the family
$\ahq$ can be considered as a q-analogue 
of commutative algebra of functions on the $\cp$ type orbits and therefore
it can be called {\em quantum or braided} $\cp$ type orbit?
In \cite{DG2} (following \cite{DS})
 $\tS$-commutative algebras were considered in such a role.

However, from representation theory point of view it is more reasonable
to consider as "quantum (braided) orbit of $\cp$ type" the  algebra
$\aq$ since it is the only algebra from the family $\ahq$ which cannot be 
represented
in a q-deformed Verma module $V_{\mu\om_1}$ with any $\mu$.
{}From this point of view it is a singular point like in a classical case.

One way more to define a version of a q-commutative algebra is discussed in
\cite{Do} (cf. footnote 4).

So,  there is no universal way to single out from the family $\ahq$ a 
braided analogue of a commutative algebra. All the above candidates for this 
role have
their own motivations.

Let us remark  that in $sl(2)$ case our approach leads to Podles' quantum 
sphere
(cf. \cite{P} where this algebra is also equipped with an  involution 
$*$). We do not consider here
the problem of a proper definition of involution operators (cf. \cite{DGR1}
for a discussion on this problem). We could only  emphasize  that in our 
 approach all representations of  algebras in question  are $\uq$-morphisms.
So, if we want to consider a $*$-representation theory of this algebra
 we must first introduce $*$-operator 
in the space $\End V$ in the spirit of super-theory: the classical property
$(AB)^*=B^*A^*$ will be failled. 

\end{document}